\documentclass{amsart}

\usepackage[left=3cm,top=3cm,right=3cm,bottom=3cm]{geometry}
\usepackage{amssymb, amsmath, amsthm}
\usepackage{verbatim}
\usepackage{comment}
\usepackage{graphicx}
\usepackage{enumerate}
\usepackage{color}

\DeclareGraphicsExtensions{.pdf,.jpg,.png}

\theoremstyle{plain}
\newtheorem{theorem}{Theorem}
\newtheorem{lemma}[theorem]{Lemma}

\newtheorem{corollary}[theorem]{Corollary}
\newtheorem{claim}[theorem]{Claim}
\newtheorem{conjecture}[theorem]{Conjecture}
\newtheorem{problem}[theorem]{Problem}

\theoremstyle{definition}

\DeclareMathOperator{\lpf}{lpf}

\newcommand{\eps}{\varepsilon}

\newcommand{\A}{\mathcal{A}}

\newcommand{\B}{\mathcal{B}}

\newcommand{\set}[1]{\left\{#1\right\}}
\newcommand{\size}[1]{\left|#1\right|}

\makeatletter

\begin{document}

\title{Sums, products and ratios along the edges of a graph}

\author{Noga Alon}
\address{\noindent  Schools of Mathematics and  
Computer Science, Tel Aviv University, Tel Aviv 69978,
Israel and Department of Mathematics, Princeton University,
Princeton, NJ 08544, USA} 
\email{nogaa@tau.ac.il}
\thanks{Research supported in part by an ISF grant and a GIF grant}

\author{Imre Ruzsa}
\address{\noindent Alfr\'ed R\'enyi Institute of Mathematics, 
Hungarian Academy of Sciences, POB 127, H-1364 Budapest, Hungary}
\email{ruzsa.z.imre@renyi.mta.hu}
\thanks{Research supported in part by an OTKA NK 104183 grant}
\author{J\'{o}zsef Solymosi}
\address{\noindent Department of Mathematics, 
University of British Columbia, 1984 Mathematics Road, 
Vancouver, BC, V6T 1Z2, Canada}
\email{solymosi@math.ubc.ca}
\thanks{Research supported in part by a  NSERC and an OTKA NK 104183 grant}

\date{}

\begin{abstract}
In their seminal paper Erd\H{o}s and Szemer\'edi formulated conjectures
on the size of sumset and product set of integers. The strongest form of
their conjecture is about sums and products along the edges of a graph. In
this paper we show that this strong form of the Erd\H{o}s-Szemer\'edi
conjecture does not hold. We give upper and lower bounds on the
cardinalities of sumsets, product sets and ratio sets along 
the edges of graphs.
\end{abstract}

\maketitle

\section{Introduction}

\subsection{Sum-product problems}
Given a finite set $\A$ of a ring, 
the \emph{sumset} and the \emph{product set} are defined by
$$\A+\A=\set{A+B:A,B\in\A},$$
and
$$\A\A=\set{AB:A,B\in\A}.$$

Erd\H{o}s and Szemer\'{e}di raised the following conjectrure
\begin{conjecture}
\label{c11}
\cite{ERDSZE} 
Every finite set of integers 
$\A$  having large enough cardinality, satisfies
\begin{equation}\label{conjecture1}
\max(\size{\A+\A},\size{\A\A})\geq\size{\A}^{2-\eps},
\end{equation}
where $\eps\to0$ as $|A|\to \infty$. 
\end{conjecture}
They proved that
\begin{equation}\label{sumproduct}
\max(\size{\A+\A},\size{\A\A})=\Omega(\size{\A}^{1+\delta}),
\end{equation}
for some $\delta>0$. Here and in what follows we use
the asymptotic notation
$\Omega(\cdot), O(\cdot)$ and $\Theta(\cdot)$. 
For two functions over the reals, $f(x)$ and
$g(x),$ we write  $f(x)=\Omega(g(x))$ if there is a positive constant,
$B>0,$ and a threshold, $D,$ such that $f(x)\geq B\cdot g(x)$ for all
$x\geq D.$ We write  $f(x)=O(g(x))$ if there is a positive constant,
$B>0,$ and a threshold, $D,$ such that $f(x)\leq B\cdot g(x)$ for all
$x\geq D.$ Finally, $f(x)=\Theta(g(x))$ if $f(x)=O(g(x))$ and
$f(x)=\Omega(g(x))$.

\medskip
Erd\H{o}s and Szemer\'edi formulated an even stronger conjecture. In this
variant one considers a subset of the 
possible pairs in the sumset and product
set. Let $G_n$ be a graph on $n$ vertices, $v_1, v_2, \ldots, v_n,$ with
$n^{1+c}$ edges for some real $c>0$. Let $\A$ be an $n$-element
set of real numbers, $\A=\{a_1, a_2,\ldots, a_n\}.$ The {\em sumset of
$\A$ along $G_n$}, denoted by $\A+_{G_n}\A$, is the set $\{a_i+a_j |
(i,j)\in E(G_n)\}.$ The product set along $G_n$ is defined similarly,
$$\A\cdot_{G_n}\A=\{a_i\cdot a_j | (i,j)\in E(G_n)\}.$$
The Strong Erd\H{o}s-Szemer\'edi Conjecture is the following. 

\begin{conjecture}\label{Strong}\cite{ERDSZE}
For every $c>$ and $\varepsilon >0,$ there is a threshold, $n_0,$
such that if $n\geq n_0$ then for any $n$-element subset of integers
$\A\subset \mathbb{N}$ 
and any graph $G_n$ with $n$ vertices and at least $n^{1+c}$ edges
$$
|\A+_{G_n}\A|+|\A\cdot_{G_n}\A|\geq |\A|^{1+c-\varepsilon}.
$$
\end{conjecture}

The original conjecture, inequality (\ref{conjecture1}), would follow
from this stronger conjecture by taking the complete graph, $G_n=K_n.$

For more details on the sum-product problem, we refer to a recent survey
\cite{GRASOL}.

Here we refute Conjecture \ref{Strong} by giving constructions with
small sumsets along a graph where the product set is also small. A similar
problem -- which is closely related to the original sum-product conjecture
-- is to bound the number of sums and {\em ratios} along the edges of
a graph. We give upper and lower bounds on these quantities.

\section{Products}
In the next construction we define a set and a graph with 
many edges such that both the sumset and the product set are small.

\begin{theorem}
\label{sumprod}
For arbitrary large $m_0,$ there is a set of integers, $\A,$ and graph
on $|\A|=m\geq m_0$ vertices, $G_m,$ with $\Omega(m^{5/3}/\log^{1/3}m)$
edges such that
$$
|\A+_{G_m}\A|+|\A\cdot_{G_m}\A|= O\left((|\A|\log |\A|)^{4/3}\right).
$$
\end{theorem}

{\em Proof:} \,
It is easier to describe our construction using rational numbers instead
of integers. Multiplying then with the least common multiple of the
denominators will not effect the size of the sumset or the product set,
giving a construction for integers.

We define the set $\A$ 
first and then the graph. Below, the function $\lpf(m)$
denotes the least prime factor of $m.$ We write $(v,w)=1$ if $v$ and $w$
are relatively prime.

$$ 
\A:=\left\{\frac{uw}{v} \text { } | \text { } u,v,w\in\mathbb{N},
\text{ where } v,w\leq n^{1/6}, u\leq n^{2/3}, (v,w)=1,\text{ and }
\lpf(u)> n^{1/6} \right\}. 
$$

The number of $u,w,v$ triples with $v,w\leq n^{1/6}, u\leq n^{2/3}$
is about $n,$ but there are further restrictions. The $\lpf(u)>n^{1/6},
u\leq n^{2/3}$ conditions allow us to select about 
$n^{2/3}/\log (n^{1/6})$
numbers for $u,$ and there are $\sim 6n^{1/3}/\pi^2$ coprime $v,w$
pairs up to $n^{1/6}.$ We are going to define a graph $G_{m}$ with
vertex set $\A,$ where $|\A|=m= O(n/\log n)$. Two elements, $a,b\in\A$
are connected by an edge if in the definition of $\A$ above $a=\frac{wu}{v}$
and $b=\frac{vz}{w}.$ There are at least 
$$
\Omega(n^{1/6}n^{1/6}(n^{2/3}/\log
(n^{1/6}))^2)=\Omega(n^{5/3}/\log^2n)
$$ 
edges (the number of quadruples $u,v,w,z$
satisfying $v,w\leq n^{1/6}, u,z\leq n^{2/3}, (v,w)=1,$ and $\lpf(z),
\lpf(u)> n^{1/6}$).

The products of pairs of elements of $\A$ along an edge of
$G_{m}$ are integers of size at most $n^{4/3}.$ The sums
along the edges are  of the form
$$\frac{wu}{v}+\frac{vz}{w}=\frac{w^2u+v^2z}{vw}.$$
The denominator is a positive integer of size at most $n^{1/3}$ and
the numerator is a positive integer of size at most $2n$, hence
the number of sums is at most
$2n^{4/3}.$ \qed

\medskip
Modifying the construction above we give a counterexample to the
Strong Erd\H{o}s-Szemer\'edi Conjecture for every $1>c>0.$ For the
sake of simplicity we will ignore logarithmic multipliers, using the
asymptotic notations $\Omega_l(\cdot), O_l(\cdot)$ and $\Theta_l(\cdot)$. 
For two functions
over the reals, $f(x)$ and $g(x),$ we write  $f(x)=\Omega_l (g(x))$
if there is a constant, $B\leq 0,$ and a threshold, $D,$ such that
$f(x)\geq \log^Bx\cdot g(x)$ for all $x\geq D.$ We write  $f(x)=O_l
(g(x))$ if there is a positive constant, $B\geq 0,$ and a threshold,
$D,$ such that $f(x)\leq \log^Bx\cdot g(x)$ for all $x\geq D.$
$f(x)=\Theta_l (g(x))$ if $f(x)=O_l(g(x))$ and
$f(x)=\Omega_l(g(x)).$

\begin{theorem}
For every $1>c>0$  there is a $\delta >0,$ such that for arbitrary
large $n$ there is an $n$-element subset of integers, $A\subset
\mathbb{N},$ and a graph $H_n$ with $\Omega_l(n^{1+c})$ edges such that
$$|\A+_{H_n}\A|+|\A\cdot_{H_n}\A|=O_l(|\A|^{1+c-\delta}).$$
\end{theorem}
{\em Proof:} \,
We consider two cases separately, when $0<c \leq 2/3$ and when $2/3<c<1.$
\begin{itemize}
    \item[]{\bf Case 1.} ($2/3<c<1$) 
    We define $\A$ similar to the previous construction, 
but now the ranges of $u,v$ and $w$ are different. 
 
$$
\A:=\left\{\frac{uw}{v} \text { } | \text { } v,w
\leq n^{\frac{1-c}{2}},(v,w)=1, u\leq n^{c}, 
\lpf(u)> n^{\frac{1-c}{2}} \right\}. 
$$

The number of $u,w,v$ triples satisfying the conditions is $\Omega_l(n).$
If $|\A|=m$ then let $H_{m}$ be the graphs with vertex set $\A.$
Two elements, $a,b\in\A$ are connected by an edge if they can be
written as $a=\frac{uw}{v}$ and $b=\frac{vz}{w}.$  There are at
least $\Omega_l(n^{1+c})$ edges in $H_{m}.$ The products of two such
elements of $\A$ are integers of size at most $n^{2c}.$ A typical sum is
$$
\frac{uw}{v}+\frac{vz}{w}=\frac{w^2u+v^2z}{vw}.
$$ 
The numerator is an integer of size at most  $2n$ and the
denominator is an integer of size at most $n^{1-c}$. 
Therefore
the sumset along the edges of $H_{m}$ has size at most $2n^{2-c}.$
Since $2c>2-c$ in this range of $c,$ we set $\delta=1-c.$ (Note that
$n=O_l(m)$)
    \item[]{\bf Case 2.} ($0<c \leq 2/3$) It is possible to describe a
construction similar to that in the first case, but we prefer to
take a subgraph of the
    graph $G_m$ in Theorem \ref{sumprod}.  Let $p$ be a parameter
satisfying  $0 <p \leq 1$, to be specified later.
    In $G_m$ take first the edges with the
    $pm^{4/3}$ most popular products. This gives a graph $G_m$ with at
    least $\Omega_l(pm^{5/3}$) edges and with $O_l(pm^{4/3})$ products
    and at most $O_l(m^{4/3})$ sums. Now in $G_m$ take the most popular
    $pm^{4/3}$ sums to get the subgraph $H_m$ with at least $\Omega_l(p^2
    m^{5/3})$ edges, with $O_l(pm^{4/3})$ products and  $O_l(pm^{4/3})$
    sums. Choosing $p$ to be $n^{c/2}/n^{1/3}$ we get $\Omega_l(n^{1+c})$
    edges and $O_l(n^{1+c/2})$ sums and products. \qed
\end{itemize}
\medskip

\section{Ratios}
In this section, we consider a problem similar to the Strong
Erd\H{o}s-Szemer\'edi Conjecture, but we change products to ratios. 
Define 
$$\A/_{G_n}\A=\{a_i/ a_j | (i,j)\in E(G_n)\}.
$$ 
(Note that each edge $(i,j)$ here providse two ratios:
$a_i/a_j$ and $a_j/a_i$. )
Changing product to
ratio is a common technique in sum-product bounds. When one is using the
multiplicative energy (like in \cite{SOL} and \cite{KONSHK} for example)
then the role of product and ratio are interchangeable. The multiplicative
energy of a set $\A$ is the number of quadruples $(a,b,c,d)\in \A^4$ such
that $ab=cd,$ which is clearly the same as the number of quadruples where
$a/c=d/b.$ But the symmetry fails in the Strong Erd\H{o}s-Szemer\'edi
Conjecture. We are going to show examples when the sumsets and ratiosets
are even smaller than in the previous construction.

\subsection{Connection to the original conjecture}
What is the connection of the Strong Erd\H{o}s-Szemer\'edi Conjecture
to the original conjecture (when $G_n=K_n$)? Similar questions were
investigated in \cite{CHM}. Here we consider the connections to the
sum-ratio problem along a graph. If there was a counterexample to
Conjecture \ref{c11}, that would imply the existence of a set
with very small sumset and 
ratioset along a dense graph. In our first result let us 
suppose that both the product set and ratioset are small. 

\begin{theorem}\label{sum-prod_fail} Let us suppose that there is a
set of $n$ real numbers, $\A,$ such that $|\A+\A|\leq
n^{2-\alpha}$,$|\A\A| =\Theta(n^{2-\beta})$
and
$|\A/\A|\leq n^{2-\beta}$  
for some $\alpha,\beta>0$ real numbers. Then
there is a set $\B$ with $N>n$ elements and a graph $G_N$ with
$\Omega(N^{\frac{3}{3-\beta}})$ edges such that   $$|\B/_{G_N}\B|=O(N),$$
and $$|\B+_{G_N}\B|=O(N^{\frac{2-\alpha}{3-\beta}}).$$ \end{theorem}

{\em Proof:} Let $\B=\{a\pm \zeta 
bc\text{ }|\text{ } a,b,c\in\A\},$ where $\zeta\in R$ is selected 
such that all sums are distinct, 
$\B$ has cardinality $2|\A||\A\A|=N=\Theta(n^{3-\beta}).$ 
In the graph, $G_N,$ 
every $a-\zeta ac$ is connected to $b+\zeta ac$ by an edge. 
The number of edges is $n^3=\Omega(N^{\frac{3}{3-\beta}}).$ 
The number of sums along the edges 
is $|\A+\A|=O(N^{\frac{2-\alpha}{3-\beta}}),$  and ratios along
the edges have the form 
\[ \frac{a+\zeta ac}{b-\zeta ac}= \frac{1+\zeta c}{b/a-\zeta c}\]
so the cardinality of the ratioset is at most $|\A||\A/\A|=O(N).$ 
\qed

A bound on the cardinality of the product set does not imply a similar
bound on the ratio set. Or, equivalently, a  bound on the cardinality of the sumset does not imply a similar
bound on the difference set. A classical construction of the second author \cite{RUZSA}
is an example for that. It uses the observation that $S = \{0, 1, 3\}$ satisfies $|S+S| = 6$ and $|S-S| = 7.$
If we consider the set of numbers, $A,$  of the form $a=\sum_{i=0}^{k-1}\alpha_i(a)10^i,$ where 
$\alpha_i(a)\in S,$ then $|A|=3^{k},$ $|A+A|=6^k,$ and $|A-A|=7^k.$ 

Note that in this construction, the multiplicity of a member
$a-a'=\sum_{i=0}^{k-1} (\alpha(a)-\alpha(a'))10^i$ of $A-A$ along the
edges of the complete graph on $A$ is $3^r$, where $r$ is the number
of indices satisfying $\alpha(a)=\alpha(a')$. It is easy to see that
for every fixed small 
$\delta>0$ the fraction of edges in which the parameter $r$ exceeds
$(1/3+\delta)k$ is at most $e^{-\Omega(\delta^2 k)}$. Therefore,
any graph on $A$ with at least 
$(9^k)^{1-c \delta^2}$ edges, for an appropriate absolute positive 
constant $c$, has on at least half of its edges a value of the difference
with multiplicity at most $3^{(1/3+\delta)k}$, implying that the number
of distinct differences along the edges is at least
$$
0.5 \frac{(9^k)^{1-c\delta^2}}{3^{(1/3+\delta)k}}.
$$
As $9/3^{1/3}  >6$ and the number of sums is only $6^k$, this shows that
for small $\delta$ the number of differences along the edges of any such
graph is
significantly larger than the number of sums.

The above discussion shows that we need a modified statement to transform a
possible counterexample to Conjecture \ref{c11} to a statement about few
sums and ratios along a graph. We are going to apply the following lemma.

\begin{lemma}\label{energy}
Let $A$ be an $n$-element subset of an abelian group and suppose that
$|A+A|\leq K|A|.$ Then there is an integer parameter, $M,$ and a graph,
$H_n,$ with vertex set $A,$ and at least $\sqrt{ \frac{M|A|^{3}}{4K\log
|A|}}$ edges, such that
$|A-_{H_n}A|\leq M.$ Moreover, $M$ satisfies the following inequalities: 
$$\frac{|A|}{K\log|A|}\leq M\leq 4K\log|A||A|.$$ 
\end{lemma}

{\em Proof:} 
The additive energy of $A,$ denoted by $E(A),$ is the number of $a,b,c,d$
quadruples from $A$ such that $a+b=c+d.$ This is the same as the number
of quadruples satisfying $a-c=d-b.$ By the Cauchy-Schwartz inequality
$E(A)\geq |A|^3/K.$ Denote the elements of the difference
set as follows
$A-A=\{t_1,t_2,\ldots,t_\ell\}.$ For every element we can define its
multiplicity, $m(t_i)=|\{a,b\in A | a-b=t_i\}|.$ With these notations
we can write the additive energy as
\[E(A)=\sum_{i=1}^\ell m^2(t_i)=\sum_{k=1}^{\log n}
\sum_{2^k\leq m(t_i)<2^{k+1}}m^2(t_i).\] 
There is a $k$ such that 
\begin{equation}\label{edge}
\sum_{2^k\leq m(t_i)<2^{k+1}}m^2(t_i)\geq \frac{|A|^3}{K\log |A|}.
\end{equation}
Let $T_k=\{t_i\in A-A \text{ }|\text{ } 2^k\leq m(t_i)<2^{k+1}\},$ 
and set $M=|T_k|.$ The edges of $H_n$ are defined as follows;
$(a,b)\in A^2$ is an edge iff $a-b=t_i\in T_k.$ The number of 
the edges is $\sum_{t_i\in T_k}m(t_i).$ From inequality (\ref{edge})
we have a lower bound on $m(t_i)$-s,
\[m(t_i)\geq \frac{|A|^{3/2}}{2\sqrt{MK\log |A|}},\]
so the number of edges is at least 
\[\sqrt{ \frac{M|A|^{3}}{4K\log |A|}}.\]
In order to bound the magnitude of $M,$ note that since 
$\sum_{t_i\in T_k}m(t_i)\leq |A|^2,$ 
the largest $m(t_i)$ for an element $t_i$ in $T_k$ satisfies the trivial
inequality $\max_{t_i\in T_k} (m(t_i))\leq 2|A|^2/M.$ 
Replacing $m(t_i)$-s by $2|A|^2/M$ 
on the left side of  inequality (\ref{edge})
we get the desired upper bound on $M.$ The lower bound follows from
the 
same inequality and from the 
fact that $m(t_i)\leq |A|$ for every $t_i\in A-A.$

\qed

\medskip
In the proof of Theorem \ref{sum-prod_fail}, in $G_N,$ the edges were
defined by pairs of vertices having the form $(a-\zeta ac, b+\zeta ac).$
If we have a bound on the product set only, $|\A\A|\leq n^{2-\beta},$
then in our new graph, $G_N',$ we connect $a-\zeta ac$ and $b+\zeta
ac$ only if $(a,b)$ is an edge in $H_n$, where $H_n$ is defined in Lemma
\ref{energy} applied to the set $\A$ in the multiplicative group.
This guarantees that the ratio set along the edges is not
(much) larger than the product set along edges of $G_N.$ The new
graph $G_N'$
is a subgraph of the graph $G_N$ in Theorem \ref{sum-prod_fail}.

The new parameter, $M,$ makes the description of our next result a bit
complicated, but the important feature of this construction is that it
shows that if there is a counterexample to Conjecture \ref{c11}
then there is a set of numbers, $\B,$ and graph, $G_N',$ with many edges
so that the sumset and the ratio set are both small.  The number of
edges might be less than in Theorem \ref{sum-prod_fail}, but then the
size of the ratio set along this graph  is much smaller. We will state
a simpler, but weaker statement in a corollary below.

\begin{theorem}\label{product_ratio}
Let us suppose that there is a
set of $n$ real numbers, $\A,$ such that $|\A+\A|\leq n^{2-\alpha}$ and
$|\A\A| = \Theta_l (n^{2-\beta})$  
for some $\alpha>0,\beta>1/2$ real numbers. Then
there is a set $\B$ with $N>n$ elements, a parameter $M$ in the range 
$$\Omega_l(N^{\frac{\beta}{3-\beta}})\leq 
M\leq O_l(N^{\frac{2-\beta}{3-\beta}}),$$ 
and a graph $G_N'$ with
$\Omega_l(M^{\frac{1}{2}}N^{\frac{4+\beta}{6-2\beta}})$ 
edges such that   $$|\B/_{G_N'}\B|=O_l(MN^{\frac{1}{3-\beta}}),$$
and $$|\B+_{G_N'}\B|=O_l(N^{\frac{2-\alpha}{3-\beta}}).$$
\end{theorem}

Note that the number of edges in $G_N'$ is at least
$\Omega_l(N^{\frac{2+\beta}{3-\beta}})$ which is bigger than the
cardinalities of the sumset and ratio set along its edges.

\begin{corollary}
Let us suppose that there is a
set of $n$ real numbers, $\A,$ such that $|\A+\A|\leq n^{2-\alpha}$ and
$|\A\A| =\Theta_l(n^{2-\beta})$  
for some $\alpha>0,\beta>1/2$ real numbers. Then
there is a set $\B$ with $N>n$ elements, and a graph $G_N'$ with
$\Omega_l(N^{\frac{2+\beta}{3-\beta}})$ 
edges such that   $$|\B/_{G_N'}\B|=O_l(N),$$
and $$|\B+_{G_N'}\B|=O_l(N^{\frac{2-\alpha}{3-\beta}}).$$
\end{corollary}

\medskip
\noindent
{\em Proof of Theorem \ref{product_ratio}:}  Applying
Lemma \ref{energy} to the multiplicative subgroup
of real numbers with $K=|\A|^{1-\beta}$ we get a graph $H_n$ as
in the lemma.
We connect $a-\zeta ac$
and $b+\zeta ac$ in $G_N'$ only if $(a,b)$ is an edge in $H_n.$ For
given $a$ and $b$ we can choose any $c\in\A,$ so the number of edges is
$\Omega_l(|\A|\sqrt{M|\A|^{2+\beta}})=\Omega_l(\sqrt{M|\A|^{4+\beta}}).$
The size of the ratio set along the edges is $|\A|M,$ and the sumset
is not larger than $O(N^{\frac{2-\alpha}{3-\beta}})$ since $G_N'$ is a
subgraph of $G_N.$
\qed

\subsection{Constructions}

In the next construction we define a set and a graph with many edges
such that both the sumset and the ratio set are very small. It can be
viewed as a special case of Theorem \ref{sum-prod_fail} with $\alpha=0$
and $\beta=1.$

\begin{theorem}\label{projection}
For arbitrary large $n,$ there is a set of reals 
and graph, $G_n$ with $\Omega(n^{3/2})$ edges such that
$$|\A+_{G_n}\A|+|\A/_{G_n}\A|\leq O(|\A|).$$
\end{theorem}
{\em Proof:} \,
As before, in the construction we define the set $\A$ first and then
$G_n.$ Let $\A=\{\pm(2^i-2^j) | 1\leq j<i\leq \sqrt{n}\}.$ Two elements,
$2^i-2^j$ and $-(2^k-2^\ell),$ are connected by an edge iff $j=\ell.$
Along this $G_n$ both the sumsets and the ratiosets are small,
$$|\A+_{G_n}\A|=|\{2^i-2^k | 1\leq i,k\leq \sqrt{n}\}|\leq n,$$
and
$$|\A/_{G_n}\A|=|\{-(2^s-1)/(2^t-1) | 1\leq s,t\leq \sqrt{n}-1\}|<n.$$
In this construction $$|\A|=2{\lfloor\sqrt{n}\rfloor\choose
2}\sim n,$$ and the number of edges is a little more than
$$2{\lfloor\sqrt{n}\rfloor\choose 3}
\geq \left(\frac{1}{3}-o(1)\right)n^{3/2}.$$ \qed
\medskip 

\subsection{Matchings} 
Erd\H{o}s and Szemer\'edi mentioned in their paper that maybe even for
a linear number of edges (when $c=0$) the Strong Erd\H{o}s-Szemer\'edi
Conjecture holds, but noted that it is not true for reals. There
are sets of reals such that $G_n$ is a perfect matching and
$$|\A+_{G_n}\A|+|\A\cdot_{G_n}\A|= O(|\A|^{1/2}).$$
It was shown by Alon,  Angel, Benjamini, and Lubetzky in  \cite{AABL}
that if we assume the Bombieri-Lang conjecture (see details in
\cite{CHM}), then for any set of integers $\A$, if $G_n$ is a
matching, $$|\A+_{G_n}\A|+|\A\cdot_{G_n}\A|=
\Omega(|\A|^{4/7}).$$ It is possible that $4/7$ can be improved 
to a number close
 to 1, but if we change multiplication to ratio, then just the trivial
 bound, $\Omega(\sqrt{n}),$ holds.

A simple construction demonstrating 
this is the following.  Take $n=k^2$ and distinct
primes, $p_1,.,p_k, q_1,..,q_k$. The matching consists of all pairs
$(p_i/q_j, (q_j-1)p_i/q_j)$   $(i,j=1,..,k).$  Then all these rationals
are pairwise distinct, the sums along the matching edges are the $p_i$ s,
the quotients (of large divided by small) along the edges are $(q_j-1).$
For this set $\A$ and matching $G_n$ we have 
$$|\A+_{G_n}\A|+|\A/_{G_n}\A|=
O(|\A|^{1/2}),$$ 
which is as small as possible.

\section{Lower bounds} 
Lower bounds on the number of sums and products along graphs were
obtained in \cite{AABL}. Under assuming the Bombieri-Lang conjecture
they proved  that if $\A$ is an  $n$-element set of integers and $G_n$
a graph with $m$ edges then
\begin{equation}\label{Bomb}
    |\A+_{G_n}\A|+|\A\cdot_{G_n}\A|= \Omega\left(
\min\left(\frac{m^{8/14}}{n^{1/14}},\frac{m}{n^{1/2}}\right)\right).
\end{equation}

For the unconditional 
case (without the Bombieri-Lang conjecture) they proved that
\begin{equation}\label{uncond}
    |\A+_{G_n}\A|+|\A\cdot_{G_n}\A|= 
\Omega\left(\frac{m^{19/9-o(1)}}{n^{28/9+o(1)}}\right).
\end{equation}
We will apply a variant of Elekes'  proof,
used in his sum-product estimate in \cite{ELE} to get better estimates. 

\begin{theorem}\label{bound_product}
Let $\A$ be an  $n$-element 
set of reals and $G_n$ a graph with $m$ edges. Then
$$|\A+_{G_n}\A|+|\A\cdot_{G_n}\A|\geq 
\Omega\left(\frac{m^{3/2}}{n^{7/4}}\right).$$
\end{theorem}
{\em Proof:} Let us consider the Cartesian product,
$(\A+_{G_n}\A)\times(\A\cdot_{G_n}\A),$ where the sums and the
products are considered along a graph $G_n.$ Define a set of
$n^2$ lines $L,$ where the lines are $y=(x-a)b$ for every $a,b \in \A.$
For any $u,$ an element of $\A,$ if in the graph $u$ has two neighbours
$w_1,w_2,$ then $((u+w_1)-w_1)w_2$ is in $ \A\cdot_{G_n}\A.$ Thus
$(u+w_1,uw_2)$ lies on the line $y=(x-w_1)w_2$ and hence
the lines
give at least the sum of the squares of degrees incidences in the
Cartesian product $(\A+_{G_n}\A)\times(\A\cdot_{G_n}\A)$. If the
graph has $m$ edges, then by the Cauchy-Schwartz inequality the
number of incidences is at least $n(2m/n)^2.$ On the other hand,
by the Szemer\'edi-Trotter theorem \cite{SzeTr}, we have
at most $O(n^{4/3}(|\A+\A||\A/\A|)^{2/3})$ incidences. We conclude
that
$m^2/n<
cn^{4/3}(|\A+\A||\A/\A|)^{2/3}.$ This implies  the 
required $\Omega((m^6/n^7)^{1/4})$
lower bound on $|\A+_{G_n}\A|+|\A\cdot_{G_n}\A|.$
\qed

\medskip
Since the values of the product and  sum along an edge determine
the values in the end-vertices, there is an obvious lower bound,
$|\A+_{G_n}\A|+|\A\cdot_{G_n}\A|\geq \sqrt{m}.$ Note that Theorem
\ref{bound_product} gives stronger bound only if the number of edges is
larger than $n^{7/4}.$ Our result improves the (conditional) inequality
in (\ref{Bomb}) if the number of edges is larger than $n^{47/26}\sim
n^{1.8},$ and it is always stronger than the bound in (\ref{uncond}).

\medskip
The very same technique can be applied to give a similar lower bound
on the number of sums and ratios. Consider the Cartesian product,
$(\A+_{G_n}\A)\times(\A/_{G_n}\A),$ and the set of lines $L,$ where
the lines are $y=(x-a)/b$ for every $a,b \in \A.$ Applying the
Szemer\'edi-Trotter theorem as above, we have

$$|\A+_{G_n}\A|+|\A/_{G_n}\A|
\geq \Omega\left(\frac{m^{3/2}}{n^{7/4}}\right).$$

\medskip
Elekes' bound was improved in \cite{SOLY_2}, using the Szemer\'edi-Trotter
theorem in a different way. The argument there can be modified to
bound the number of sums and ratios along a graph. The proof is rather
technical, it follows \cite{SOLY_2} step by step, but with more parameters
in order to deal with the density version of the original proof. We
do not
think that this estimate is close to the truth and it is just slightly
better, in a small range when $m\gg n^{11/6},$ than the simple bound
above. We state the bound without the detailed proof.

\begin{claim}\label{bound_ratio}
Let $\A$ be an  
$n$-element set of reals and $G_n$ a graph with $m$ edges. Then
$$|\A+_{G_n}\A|+|\A/_{G_n}\A|
\geq \Omega\left(\frac{m^{18/11}}{n^{2}}\right).$$
\end{claim}

\section{Arrangements of pencils}
The following question was asked by Misha Rudnev \cite{RUD}. An
$n$-pencil in the plane is a set of $n$ concurrent lines. The center of
a pencil is the common intersection point of its lines.

\begin{problem}
If the centers of four $n$-pencils are not collinear, then what is the 
maximum possible number of points with four incident 
lines (one from each pencil)?  
\end{problem} 

Chang and Solymosi proved in \cite{CHASOL} that the number of such points
is at most $O(n^{2-\delta})$ for some $\delta>0.$ (They did not calculate
$\delta$ explicitly.) Using the construction in Theorem \ref{projection}
we show that $\delta\leq 1/2.$


\begin{claim}
For arbitrary large $n,$ there are arrangements of four non-collinear 
$n$-pencils which determine $\Omega(n^{3/2})$ points incident to four lines.
\end{claim}

{\em Proof:} 
In this construction we refer to Theorem \ref{projection}. If a set
of reals, $A,$ has small sumset then the geometric interpretation of
this fact is that the points of the Cartesian product $A\times A$ can
be covered by a small number of slope $-1$ lines. Similarly, if the
ratio-set is small then $A\times A$ can be covered by a small number
of lines through the origin.
The set of points where four lines intersect is defined as
\[P:=\{(2^i-2^j,-(2^k-2^j))
\in \mathbb{R}^2 | 1\leq j\leq i,k\leq \sqrt{n}  \}.\]
The four pencils are 
\begin{itemize}
    \item[] The vertical lines with a point in $P,$
    \[L_1:=\{x=2^i-2^j | 1\leq j\leq i\leq \sqrt{n}  \}.\]
    \item[] The horizontal lines with a point in $P,$
    \[L_2:=\{y=-2^i+2^j | 1\leq j\leq i\leq \sqrt{n}  \}.\]
    \item[] The slope $-1$ lines with a point in $P,$
    \[L_3:=\{x-(2^i-2^j)=-(y+(2^k-2^j)) | 1\leq j\leq i,k\leq \sqrt{n}  \},\]
    \item[] Lines through the origin with a point in $P,$
    \[L_4:=\left\{y=-\frac{2^{i-j}-1}{2^{k-j}-1}x | 
1\leq j\leq i,k\leq \sqrt{n}  \right\}.\]
\end{itemize}
Note that in the definition of $L_3$ and $L_4$ the same lines are
listed multiple times. Ignoring these repetitions it is easy to see
that all four families have size approximately $n,$
and $|P|=\Omega(n^{3/2}).$ One can apply a projective transformation to
shift the centers of the pencils from infinity to $\mathbb{R}^2.$ \qed

\section*{Acknowledgement}
We are thankful to Misha Rudnev, whose question on pencils initiated
this work.

\end{document}